\let\documentclass\relax\documentclass % Для тупого арХива
\input AHTOHFIE.STY
\hfuzz6.8pt
\def\r{{\bf\=r}}
\def\rk{\={\rm rk}}    %переопределён
\def\rank{{\rm rank}} %переопределён слегка

\UDC{
%512.543.12%  %Свободные группы
512.543.52% %Свободные произведения
+512.543.14% %Группы с конечным числом образующих
+512.544.23% %Структурные свойства групп
}


\title{%
An analogue of the strengthened Hanna Neumann conjecture
for virtually free groups\\
and virtually free products
%05.06.21
%of left-orderable groups
%%
%Пересечение подгрупп в почти свободных произведениях
%левоупорядочиваемых групп
}
\author{%
Anton A. Klyachko$^\flat$
\quad
Alexander O. Zakharov$^\sharp$
}
\address{
$^\flat$%
%\мойАдресШ
Faculty of Mechanics and Mathematics of Moscow State University
\\
Moscow 119991, Leninskie gory, MSU.
\\
Moscow Center for Fundamental and Applied mathematics.
\\
klyachko@mech.math.msu.su
\\
$^\sharp$%
Instytut Matematyczny,
Uniwersytet Wroc\l awski,
pl. Grunwaldzki 2/4,
50-384, Wroc\l aw.
%POLAND
%Исследовательская лаборатория им. П. Л. Чебышева,
%Санкт-Петербург 199178, 14 линия В.О., дом 29Б
\\
alexander.zakharov@uwr.edu.pl
}
 
\grantsFirst{%
%28.05.23
%\RFBR19-01-00591,
%and the Ministry of Education and Science of the Russian Federation as 
%part of the program of the Moscow Center for Fundamental and Applied 
%Mathematics under the agreement 075-15-2019-1621 
\RSF 22-11-00075%
}

\grantsSecond{(Polish) Narodowe Centrum Nauki, grant UMO-2018/31/G/ST1/02681}

\abstract{%
The Friedman--Mineyev theorem, earlier known as the
(strengthened)
Hanna Neumann conjecture,
%%говорит, что
%%$\rank(A\cap B)-1\le(\rank(A)-1)(\rank(B)-1)$
%%для любых нетривиальных подгрупп $A$ и $B$ свободной группы.
gives a sharp estimate for the rank of the intersection of two subgroups
in a free group. We obtain an
analogue
of this inequality for any two
subgroups in a
%$A$ и $B$
virtually free group (or, more generally, in a group containing a free
product of left-orderable groups as a finite-index subgroup).
%и в почти свободных группах.
}

%%%%%%%%%%%%%%%%%%%%%%%%%
\s 0.
Introduction

The Hanna Neumann Conjecture (1957), proved
independently by
Mineyev
%([Mi12a], [Mi12b])
and Friedman
%[Fr14],
is the following fact.

\proclaim Friedman--Mineyev theorem{\rm[Mi12a], [Mi12b], [Fr14]}.
For any nontrivial subgroups $A$
and $B$
of a free group $F$,
$$
\rank(A\cap B)-1\le (\rank(A)-1)\cdot(\rank(B)-1);
\eqno{\hbox{\emph{\(the classical Hanna Neumann conjecture\)}}}
$$
moreover, for any system of representatives $S$ of the
double cosets~$AsB$ in~$F$,
$$
\sum\limits_{s \in S}\=\rank(A\cap sBs^{-1})
\le
\=\rank(A)\cdot\=\rank(B),
\eqno{\hbox{\emph{\(the strengthened Hanna Neumann conjecture\)}}}
$$
where $\=\rank(H)\:=\max(0,\;\rank(H)-1)$ is the \emph{reduced rank}
of a
free group $H$.

Alternative proofs
and generalisations of this result can be found,
e.g., in
[D12],
[AMS14],
[Za14],
[ASS15],
[Nos16],
[HW16],
[Iv17],
[JZ17]
and
[KP20].
In particular, the following analogue of the
classical Hanna Neumann conjecture for
free subgroups of a virtually free group was obtained in [KP20]:
\disp{\sl
for any free subgroups $A$ and $B$ of a virtually free group~$G$
containing a free finite-index subgroup $F$
\newline\centerline{}
\centerline{\strut
$
\=\rank(A\cap B)
\le
|G{:}F|\cdot\=\rank(A)\cdot\=\rank(B).
$
}
}%
This estimate strengthens earlier known inequalities [Za14], [ASS15]
(and is already sharp). We generalise this fact in two directions:
\-
first, we obtain an analogue of the \spacing{strengthened} Hanna Neumann
conjecture;
\-
and secondly, our estimate applies to \spacing{arbitrary} subgroups
$A$ and $B$ of a virtually free group.

\Theorem on intersection of subgroups in virtually free groups.
For any subgroups $A$ and $B$ of a virtually free group
$G$ containing a free group $F$ as a finite-index subgroup
and for any system of representatives $S$ of the
double cosets~$AsB$ in~$G$,
$$
\sum\limits_{s \in S}\rk(A\cap sBs^{-1})
\le
|G{:}F|\cdot
\rk(A)\cdot\rk(B).
\qbox{In particular,
$
\rk(A\cap B)
\le
|G{:}F|\cdot
\rk(A)\cdot\rk(B).
$
}
$$

\noindent
Here $\rk(H)$ is the \emph{virtual reduced rank}
of a
virtually free group\/\:
$
\rk(H)\:={1\over|H{:}K|}\cdot\max\bigl(0{ },\rank(K)-1\bigr),
$
where $K$ is a finite-index free subgroup of $H$.
It is easy to show that the virtual reduced rank is well-defined
(i.e. it does not depend on the choice of a free subgroup $K$);
and
$\rk(H)=\=\rank(H)$,
if $H$ is
free.
%06.06.21
Note that the virtual reduced rank of a virtually free
group coincides with the \emph{rank gradient}
of this group [La05].

Actually, the
theorem above
is a special case of a more general \emph{main theorem}
of
this paper (see the next section),
which is about intersections of
subgroups in virtually free products.
In particular, our main theorem generalises the
following known analogue of the
strengthened Hanna Neumann conjecture.

\Th AMS {\rm[AMS14] (see also [Iv17])}.
For any subgroups $A$ and $B$ of a free product
$G=\zvezda_{i\in I}G_i$ of left-orderable groups $G_i$
and for any system of representatives $S$ of
double cosets~$AsB$ in~$G$
$$
\sum\limits_{s \in S}\=\rank_{\rm K}(A\cap sBs^{-1})
\le
%|G{:}F|\cdot
\=\rank_{\rm K}(A)\cdot\=\rank_{\rm K}(B).
\qbox{In particular,
$
\=\rank_{\rm K}(A\cap B)
\le
%|G{:}F|\cdot
\=\rank_{\rm K}(A)\cdot\=\rank_{\rm K}(B).
$
}
$$

\noindent
Here, $\=\rank_{\rm K}(H)$ is the \emph{reduced Kurosh rank} of a
subgroup $H\subseteq G=\zvezda_{i\in I}G_i$, which is defined as follows:
the subgroup~$H$
decomposes (by the Kurosh theorem) into a free product
$H=\(\zvezda\limits_{j\in J}H_j\)*F$, where each
$H_j$ is nontrivial and conjugate to a subgroup of one of $G_i$,
and $F$ is free and trivially intersects all
conjugate to subgroups $G_i$; then
$\=\rank_{\rm K}(H)\:=\max(0,\;|J|+\rank(F)-1)$.

Our proof of the main theorem is based on Mineyev's approach [Mi12b],
but our definitions are somewhat different; therefore,
we prove everything ``from scratch"; thus,
this paper contains yet another alternative
(simpler) proof of the Friedman--Mineyev theorem.

\noindent
The main features of our argument are that,
considering actions of groups on forests,
we
\-
never refer explicitly to the quotient graph of this action
\-
and never require the cocompactness of the action.

\enditem
%12.07.21
%This allows us to say that our main theorem and all its
%corollaries stated above (including the Friedman--Mineyev theorem)
%are, in a sense, special cases of
%a quite elementary lemma on group actions
%(on sets),
%see Section~2.
%%
%12.07.21
The authors thank an anonymous referee for useful remarks.
%%
%28.05.23
The first author thanks also the Theoretical Physics and Mathematics
Advancement Foundation ``BASIS".
%%

%%%%%%%%%%%%%%%%%%%%%%%%%%%%
\s 1.
Main theorem

If a group $G$ contains a free product
$F=\zvezda_{i\in I} G_i$ of
%06.06.21
%nontrivial
infinite
groups $G_i$
as a finite-index subgroup,
then, for any subgroup $H\subseteq G$, the
\emph{virtual reduced Kurosh rank $\rk(H)$} with respect to the family
of
subgroups $G_i$ is defined as
$$
\rk(H)\:={\=\rank_{\rm K}(K)\over|H{:}K|},
$$
where $K$ is a subgroup of
finite index in $H$ contained in $F$, and
$\=\rank_{\rm K}(K)$ is the (usual) reduced Kurosh rank of a subgroup
$K$ of
$F=\zvezda_{i\in I} G_i$.
This value is well-defined, i.e. it does not depend on the choice
of $K$
%06.06.21
(because of an analogue of the Schreier formula
for the Kurosh rank [Ku83]),
but is
not conjugation-invariant, i.e. the numbers $\rk(H)$ and $\rk(gHg^{-1})$
may differ. To remedy this unpleasant feature, we define the
\emph{total virtual reduced Kurosh rank $\r(H)$}
\(with respect to the family $\{G_i\;|\;i\in I\}$\)
as:
$
\r(H)=\sum\limits_{j=1}^n\rk(g_jHg_j^{-1}),
$
where $\rk$ is the virtual reduced Kurosh rank
with respect to the given family
of subgroups, and $g_1,\dots,g_n$ are
representatives of the
%09.06.21
%left
right
cosets of $F$ in $G$.
It is easy to see that this value is conjugation-invariant and does not
depend on the choice of representatives $g_j$.
%12.07.21
Note that $\r(H)=0$ for finite $H$.

\proclaim{Main theorem}.
Let a group $G$ be a virtually free product of
left-orderable groups,
i.e. $G$ contains a finite-index subgroup
$F=\zvezda_{i\in I} G_i$, where all
groups~$G_i$ are left-orderable. Let $A$ and $B$ be
subgroups of $G$, and let $S$ be
a set of distinct
representatives of double cosets $AgB$ in $G$.
Then
$
\sum\limits_{s \in S}\r(A\cap sBs^{-1})
\le
\r(A)\cdot\r(B),
$
where $\r(H)$ is the total virtual
reduced
Kurosh rank of
$H\subseteq G$ \(with respect to the family
$\{G_i\;|\;i\in I\}$\).
\newline
In particular,
$
\r(A\cap B)
\le
%|G{:}F|\cdot
\r(A)\cdot\r(B).
$

\noindent
This generalises earlier known results:
\-
the case, where $F=G$, of our theorem is
%основной результат
Theorem AMS [AMS14]
(if, in addition,
all $G_i$ are infinite
cyclic, then we obtain the Friedman--Mineyev theorem,
earlier known, as the strengthened Hanna Neumann conjecture);
\-
in  the case, where $A$ and $B$ are free groups
trivially intersecting subgroups conjugate to free factors
$G_i$,
the assertion ``In particular" is the main result
of~[KP20].

\enditem
To derive the theorem on intersections of subgroups in virtually
free groups (see Introduction)
from the main theorem,
it suffices to note that
the virtual reduced rank $\rk(H)$
of a virtually free subgroup $H\subseteq G$
%, как он определён в этом следствии,
coincides with the virtual reduced Kurosh rank with respect to
any family of infinite cyclic subgroups whose free product
is $F$. Therefore, all terms in the definition of the total
virtual rank $\r(H)$ are equal
(and their number is the index of $F$),
i.e. $\r(H)=|G{:}F|\cdot\rk(H)$ in this case.

%%%%%%%%%%%%%%%%%%%%%%%%%
\s 2.
Actions

Strange as it may seem, we have not succeeded
in finding
the following simple lemma in the literature.

\proclaim Orbit-intersection lemma.
Let $A$ and $B$ be subgroups of a group $G$ that
acts freely on
a set
%12.07.21
$X$,
and let
$D$ be a set of distinct representatives of double cosets $AgB$.
Then
$$
\sum_{d\in D}
\hbox{\(the number of $(A^d\cap B)$-orbits\)}
\le
\hbox{\(the number of $A$-orbits\)}
\cdot
\hbox{\(the number of $B$-orbits\)}.
$$
%05.06.21
Moreover, for any $A$-invariant set
$Y\subseteq X$
and any $B$-invariant set~$Z\subseteq X$,
$$
\sum_{d\in D}\!
\Bigl(\hbox{the number of $(A^d\cap B)$-orbits
in $(d^{-1}\o Y)\cap Z$}\Bigr)
\!\le\!
\hbox{\(the number of $A$-orbits in $Y$\)}
\cdot
\hbox{\(the number of $B$-orbits in $Z$\)}.
$$

\Proof
Suppose that $G\times X\too^\o X$ is a free action,
and $X/H$ is the set of orbits of the action of a subgroup
$H\subseteq G$.
Consider the mapping
%05.06.21
$$
\Phi\:\left\{(d,U)\;|\;d\in D,\
U\in\Bigl((d^{-1}\o Y)\cap Z\Bigr)/(A^d\cap B)\right\}
\to Y/A\times Z/B,
\quad
\bigl(d,\ (A^d\cap B)\o x\bigr)\mapsto (A\o d\o x,\ B\o x).
$$
The assertion of the lemma follows
immediately from the following observations:
this mapping is
\-
well-defined,
\itemitem{}
i.e., it does not depend on the choice of
$x$ in the $(A^d\cap B)$-orbit (obviously),
\-
and injective;
\itemitem{}
indeed, $(A\o d\circ x,\;B\o x)=(A\o d'\o x',\;B\o x')$
means that
$d'\o x'\in A\o d\o x$ and $x'\in B\o x$, i.e.
$(d'B)\cap(Ad)\ne\emptyset$ (because the action is free) and,
hence, $d'=d$ (by the definition of $D$); so,
$x'\in(A^d\o x)\cap(B\o x)=(A^d\cap B)\o x$, as required.

\enditem

%%%%%%%%%%%%%%%%%%%%%%%%%
\s 3.
Actions on forests

All graphs in this paper are directed.
Let a group $G$
act
on a forest $\Gamma$ freely on edges (i.e. the stabiliser of each
edge is trivial).
A
set~$E$ of orbits of edges of~$\Gamma$
is called
\emph{maximal
essential}
if
$E$ is
an inclusion maximal
set such that
each component of the forest
$\Gamma\setminus\bigcup E$
%12.07.21
(i.e. each component of the forest
obtained from $\Gamma$ by deleting
all edges from each orbit of the orbit-set~$E$),
which is not a component of $\Gamma$,
has a nontrivial stabiliser.
The following lemma is a simple (and probably known) fact on
groups acting on trees.

%12.07.21
Note that, actually, no component of $\Gamma$
has trivial stabiliser if the number of components is finite
and the group is infinite
(but we assume neither these conditions to hold by default).

\proclaim Kurosh-rank lemma.
A group $G$
acting on a tree $\Gamma$ freely on edges
decomposes into a free product:
%\newline
$
G=F*\(\zvezda\limits_{i\in I}G_i\)
$,
where $F$ is a free group acting on $\Gamma$ freely,
and $G_i\ne\1$ are stabilisers of some vertices;
if the Kurosh rank of this decomposition is finite
\(i.e.
$\rank(F)+|I|<\infty$\),
then the cardinality of any maximal essential set~$E$
equals the reduced Kurosh rank of this decomposition:
$
|E|=\max(0,\;\rank(F)+|I|-1).
$

\noindent{\bf Sketch of a proof.}
The first assertion is well-known.
To prove the second assertion,
for any edge $e$, consider the components
$X$ and $Y$
of the forest $\Gamma\setminus(G\o e)$ connected by $e$.
The ping-pong lemma implies immediately that
$$
G=\cases{
\St(X)*\gp{g}, &if $g\o X=Y$ for some $g\in G$
(which necessarily acts freely on $\Gamma$);
\cr\cr
\St(X)*\St(Y),&if $g\o X\ne Y$ for any $g\in G$.
}
$$
An obvious induction completes the proof
%12.07.21
(as the Kurosh rank is finite).
This lemma also follows from [AMS14]
(Theorem~2.4, using the arguments of Proposition~3.4).

\medskip

\noindent
We want to
generalise this simple fact to the
case, where $\Gamma$ is a forest
consisting of finitely many trees:
$\Gamma=T_1\sqcup\dots\sqcup T_n$.
In this case, the \emph{virtual
reduced Kurosh rank} of the (action of)
the group $G$ is naturally defined:
choose in $G$ a finite-index subgroup $H$ that stabilises
a tree $T_j$ and, therefore, decomposes into a free product
$
H=F*\(\zvezda\limits_{i\in I}G_i\)
$,
where $F$ is a free group acting on $T_j$ freely,
and $G_i\ne\1$ are stabilisers of some vertices of $T_j$;
the reduced Kurosh rank of this subgroup
(with respect to the given action on $T_j$)
is
$\rk(H)\:=\max(0,\;\rank(F)+|I|-1)$,
and the virtual reduced Kurosh rank of $G$
(with respect to the given action on $\Gamma$ and given component
$T_j$ of $\Gamma$) is naturally defined as:
$\rk_j(G)\:={1\over|G:H|}\cdot\rk(H)$. It is easy to see that this
value does not depend on the choice of the subgroup~$H$
%06.06.21
(if nontrivial vertex stabilisers are infinite),
but can depend on~$j$.
%09.06.21
We call the value $\sum\limits_j\rk_j(G)$ the
\emph{total virtual reduced Kurosh rank} of this action.

\proclaim Virtual-Kurosh-rank lemma.
%06.06.21
Suppose that a group $G$ acts
%05.06.21
freely on edges
on a
forest
$\Gamma=T_1\sqcup\dots\sqcup T_n$
consisting of trees $T_j$,
the stabiliser of each $T_j$
has finite Kurosh rank,
and nontrivial vertex stabilisers are infinite.
Then
$\sum\limits_{j=1}^n\rk_j(G)=|E|$
for each maximal essential set $E$.

\Proof
Suppose that $\Gamma=\Gamma_1\sqcup\dots\sqcup\Gamma_k$
and, on each $G$-invariant forest $\Gamma_i$,
the action of $G$ is transitive on components
(i.e., for any components
$T_l,T_m\subseteq\Gamma_i$,
there exists $g\in G$ such that $g\o T_l=T_m$).
Then $E=E_1\sqcup\dots\sqcup E_k$,
where $E_i=\{G\o e\in E\;|\;G\o e\subseteq\Gamma_i\}$ is a maximal
essential set of orbits of edges of the forest $\Gamma_i$.
Therefore, it suffices to prove assertion for the case,
where the action of $G$ on $\Gamma$ is
transitive on components of $\Gamma$.

In this case, all stabilisers $H_j=\St(T_j)$ of trees $T_j$ are conjugate
and, hence, isomorphic and act on the correspondingly trees similarly.
In particular, $\rk(H_j)$ does not depend on $j$.
Moreover, $|G{:}H_j|=n$ for all $j$
(because the length of an orbit equals the index of the stabiliser).
Therefore,
$$
\sum_{j=1}^n\rk_j(G)
=
\sum_{j=1}^n{1\over|G{:}H_j|}\cdot\rk(H_j)
=
\sum_{j=1}^n{1\over n}\cdot\rk(H_1)=\rk(H_1).
$$
On the other hand, the set of $H_1$-orbits of edges
$E'=\{G\o e\cap T_1\;|\;G\o e\in E\}$
is, obviously, maximal essential
with respect to the action of $H_1$ on $T_1$.
Therefore, $|E|=|E'|=\rk(H_1)$
(the latter equality follows from the Kurosh-rank lemma).
This completes the proof.

%%%%%%%%%%%%%%%%%%%%%%%%%
\s 4.
Actions on ordered forests

We say that a graph is
\emph{ordered} if it is equipped
with a partial
order on the set of edges inducing a
linear order on the set of edges of each connected component.

\proclaim
Induced-action lemma {\rm[KP20]}.
If a group $G$ has a subgroup $F$ of a finite index~$n$,
which acts on an ordered tree~$T$ preserving the order,
then $G$
can
act
preserving the order
on an
ordered forest consisting of $n$ trees;
the stabilisers of vertices and edges under this action are
conjugate to the stabilisers of vertices and edges under the initial action
of $F$ on $T$.

% 6/06/2021

\Proof
Let $S\ni1$ be a system of representatives of the left cosets of $F$ in $G$
(i.e.  $|S|=n$). Thus, each element $g\in G$ decomposes uniquely into a
product $g={\bf s}(g){\bf f}(g)$ of an element~${\bf s}(g)\in S$ and an
element ${\bf f}(g)\in F$.

Take the ordered forest $L=\bigcup\limits_{s\in S} sT$ consisting of
$n$ copies $sT$ of the ordered tree $T$ (edges from different
copies are incomparable) and consider the usual induced
action of $G$ on $L$:
\quad $g\o st\:={\bf s}(gs)\Bigl({\bf f}(gs)\o t\Bigr)$.
Clearly, this action satisfies all requirements.
This completes the proof.

\smallskip

An edge $e$ of an ordered forest with
an order-preserving action of a group $H$
is called
\emph{important}
%05.06.21
(or \emph{$H$-important})
if it is the maximal edge on
an bi-infinite line $T(e)$ intersecting only finitely many $H$-orbits of
edges.
%12.07.21
Note that if $K\subseteq H$, then any $K$-important edge is $H$-important.

\proclaim
Important-edge lemma.
If a group~$G$ acts on an ordered
forest $T$
preserving the order and freely on edges,
then
\-
the
set $\cal E$ of orbits of
important
edges
contains a
maximal
essential
set;
\-
each finite subset
%09.06.21
${\cal E'}\subseteq\cal E$
is contained in a maximal
essential
set.
\enditem
In particular,
the
%09.06.21
total virtual reduced Kurosh rank of this action
\-
equals
$|{\cal E}|$ if $|{\cal E}|<\infty$,
\-
is infinite if ${\cal E}$ is infinite.

\Proof
The assertion ``In particular" follows from the main assertion by the
virtual-Kurosh-rank lemma. It remains to prove the main assertion.
Take a finite subset ${\cal E'}$ of $\cal E$ and
put $E=\bigcup\cal E$ and $E'=\bigcup\cal E'$
%12.07.21
(i.e. $e\in E$ if and only if $G\o e\in\cal E$;
and similarly $e\in E'$ if and only if $G\o e\in\cal E'$;
so $E$ and $E'$ are sets of edges, while
$\cal E$ and $\cal E'$ are sets of orbits of edges).
We have to establish two facts:
\item{1)}
the stabiliser $\St(K)$ of each component $K$
of the forest~$T\setminus E$ has either a fixed point or
an invariant line in $K$;

\item{2)}
but the stabiliser of each component $K$
of the forest~$T\setminus E'$ is nontrivial
if there exists an important edge
%09.06.21
$e\in E'$
in $T$ incident to a vertex
of $K$.

\enditem
The both facts are easy to prove.

\item{1)}
{%
If 1) does not hold,
then the stabiliser of a component $K$
of
$T\setminus E$
contains a rank-two free subgroup
$F(x,y)\subseteq\St(K)$ acting freely on $K$
(because each non-dihedral group
nontrivially decomposable
into a free product contains a free subgroup that
trivially intersects the free factors;
the group $G$
cannot
be dihedral, because $G$ is
torsion-free if $T$ has at least one
edge).
Let $l_x$ and $l_y$  be invariant lines in $K$ for elements
$x$ and $y$, respectively.
The intersection of these lines is a finite graph:
either an interval, a point, or the empty set
(it cannot be a ray, as is known).
Let us connect the lines $l_x$ and $l_y$ by a path $\pi$.
Let us choose finite intervals~$p_x$~and~$p_y$ such that
$l_x=\bigcup\limits_{k\in\Z}x^k\o p_x$
and
$l_y=\bigcup\limits_{k\in\Z}y^k\o p_x$; and
let us take
the maximal edges $e_x$ and $e_y$ on
the intervals~$p_x$~and~$p_y$.
Without loss of generality, we can assume
that
\itemitem{-}
$x\o e_x<e_x$ and $y\o e_y<e_y$
(replace $x$ with $x^{-1}$ and/or $y$ on $y^{-1}$, if this is not the
case);
\itemitem{-}
$\(\bigcup\limits_{k=0}^\infty x^k\o p_x\)\cap(l_y\cup\pi)=
\emptyset=
\(\bigcup\limits_{k=0}^\infty y^k\o p_y\)\cap(l_x\cup\pi)$
(replace $p_x$ with $x^n\o p_x$ and/or $p_y$ with
$y^n\o p_y$ for sufficiently large $n\in\N$, if this not the case).

%\enditem
Let us connect now $p_x$ and $p_y$ by a path $p\supset(p_x\cup p_y)$.
The maximal edge $e$ of
$p$ the maximal edge on the
line~$
p\;
\cup
\(\bigcup\limits_{k=0}^\infty x^k\o p_x\)
\cup
\(\bigcup\limits_{k=0}^\infty y^k\o p_y\).
$
Thus the edge $e$ is important (Fig.~1).
This contradiction completes the proof of~1).
}

\goodbreak
%\vskip1cm plus 1cm minus5mm
\bigskip
\centerline{\input 1.PIC}
\nobreak%
%\vskip5mm%
\centerline{Fig. 1}%
%\vskip1cm plus 1cm minus 5mm%
\goodbreak
\bigskip
%\par
%\noindent

\item{2)}
Suppose that an important edge $e\in E'$
ends at a vertex of $K$.
If the edge~$g\o e$ ends also at
a vertex of $K$, then $g\in\St(K)$
and, therefore, $\St(K)\ne\1$ if $g\ne1$.
Hence, it suffices to consider the case,
where there are only finitely many (at most $2|\cal E'|$)
important edges from $E'$ incident to vertices of $K$.
Let $e\in E'$ be the minimal edge of $E'$ incident to a vertex of $K$.
Then an infinite ray of $T(e)$
(from the definition of the importance) must lie in $K$ (because
of the minimality of $e$). Since this ray can intersect only
finitely many orbits of edges (by the definition of importance),
we obtain an infinite set of edges of $K$ lying in the same orbit.
Thus, $\St(K)\ne\1$ as required.

%%%%%%%%%%%%%%%%%%
\s 5.
Proof of the main theorem

Put $n=|G{:}F|$ and let $T$ be the (Bass--Serre) tree for the decomposition
$F=\zvezda\limits_{i\in I} G_i$, i.e.
$F$ acts on $T$ freely on edges and in such a way that
the stabiliser of each vertex is conjugate to one of
the factors $G_i$.
The tree $T$ can be ordered:
the order on the set of edges of $T$ is induced by a
left-invariant order on group $F$
(which, as is known, exists [Vi49], [D\v S20]).
Thus, the action of $F$ on
$T$ preserves the order and is free on edges.
By the induced-action lemma, the group $G$
acts
on an ordered forest
$\Gamma=T_1\sqcup\dots\sqcup T_n$ consisting of
$n$ trees $T_j$
transitively on components,
freely on edges, and preserving the order.
Moreover, $\St(T_1)=F$ and $T_j=g_jT_1$, where $g_1=1,g_2,\dots,g_n$ are
representatives of the
left
cosets of $F$ in $G$.

The groups $A$ and $B$ act on the forest $\Gamma$ freely on edges and
preserving the order. Then
%12.07.21
$$
\eqalign{
&\sum_{s\in S}
\;\bigl(\hbox{the number of $(A^s\cap B)$-orbits
$(A^s\cap B)$-important edges}\bigr)
\le
\cr
\le
&\sum_{s\in S}
\;\bigl(\hbox{the number of $(A^s\cap B)$-orbits
of edges that are both $A^s$-important and $B$-important}\bigr)
=
\cr
=
&\sum_{s\in S}
\;\Bigl(\hbox{the number of $(A^s\cap B)$-orbits in the set
$\(s^{-1}\o\{A\hbox{-important edges}\}\)\cap\{B\hbox{-important edges}\}$}
\Bigr)
\le
\cr
\le&
\bigl(\hbox{the number of $A$-orbits of
$A$-important edges}\bigr)
\cdot
\bigl(\hbox{the number of $B$-orbits of
$B$-important edges}\bigr),
}
\eqno{(*)}
$$
where
\-
the first inequality holds because
any $H$-important edge is $G$-important if $H\subseteq G$;
\-
the equality holds because
an edge $e$ is
$A$-important if and only if
$s^{-1}\o e$ is $A^s$-important;
\-
the last inequality is
the orbit-intersection
lemma
%05.06.21
applied to
$$
Y=\{\hbox{$A$-important edges of $\Gamma$}\}
\subseteq
X=\{\hbox{edges of $\Gamma$}\}
\supseteq
Z=\{\hbox{$B$-important edges of $\Gamma$}\}.
$$

\enditem
By the important-edge lemma,
the set of orbits of important edges is
maximal essential (if a maximal essential
set is finite).
%12.07.21
Thus, the number of $C$-orbits of $C$-important edges
in $(*)$ is the cardinality of the
maximal essential set with respect to the action of the group $C$
on $\Gamma$ (where $C$ is $A$, $B$, or $A^s\cap B$).
Therefore, by the virtual-Kurosh-rank lemma
$$
\sum\limits_{s \in S}\r(A\cap sBs^{-1})
\le
\r(A)\cdot\r(B),
\qbox{where } \r(H)=\sum_{j=1}^n\rk_j(H)
$$
and $\rk_j(H)$ is the virtual reduced Kurosh rank
with respect to (the corresponding decomposition of)
$\St(T_j)$.
It remains to note that
the ``corresponding decomposition" of the
stabiliser of the $j$th tree has form
$
\St(T_j)=g_jFg_j^{-1}
=
\zvezda\limits_{i\in I}g_jG_ig_j^{-1}.
$
This completes the proof.

%%%%%%%%%%%%%%%%%%%%%%%%%
%\baselineskip 10 pt

\References

[AMS14]
Y. Antol\'\i n, A. Martino, and I. Schwabrow,
Kurosh rank of intersections of subgroups of free products of
right-orderable groups,
Mathematical Research Letters, 21:4 (2014), 649-661.
\arXiv 1109.0233

[ASS15]
V. Ara\'ujo, P. V. Silva, and M. Sykiotis,
Finiteness results for subgroups of finite extensions,
J. Algebra, 423 (2015), 592-614.
\arXiv 1402.0401

[D12]
W. Dicks,
Simplified Mineyev,
{\def~{\char`~}
\tt https://mat.uab.cat/~dicks/pub.html}\thinspace.

[D\v S20]
W. Dicks, Z. \v Suni\'c,
Orders on trees and free products of left-ordered groups,
Canadian Mathematical Bulletin, 63:2 (2020), 335-347.
%doi:10.4153/S0008439519000389
\arXiv 1405.1676

[Fr14]
J. Friedman,
Sheaves on graphs, their homological invariants, and a proof of the
Hanna Neumann conjecture. With an appendix by Warren Dicks,
Mem. Amer. Math. Soc. 233:1100 (2014).
\arXiv 1105.0129

[HW16]
J. Helfer and D. T. Wise,
Counting cycles in labeled graphs:
the nonpositive immersion property for one-relator groups,
International Mathematics Research Notices 2016:9 (2016), 2813-2827.

[Iv17]
S. V. Ivanov,
Intersecting free subgroups in free products of left ordered groups,
Journal of Group Theory, 20:4 (2017), 807-821.
\arXiv 1607.03010

[JZ17]
A. Jaikin-Zapirain,
Approximation by subgroups of finite index and the Hanna Neumann conjecture,
Duke Mathematical Journal, 166:10 (2017), 1955-1987.

%6/06/2021

[Ku83]
R. S. Kulkarni,
An extension of a theorem of Kurosh and applications to Fuchsian groups,
Michigan Mathematical Journal, 30:3 (1983), 259-272.

[La05]
M. Lackenby,
Expanders, rank and graphs of groups,
Israel Journal of Mathematics, 146:1 (2005), 357-370.
\arXiv math/0403127

[KP20]
A. A. Klyachko, A. N. Ponfilenko,
Intersections of subgroups in virtually free groups
and virtually free products,
Bull. Austral. Math. Soc., 101:2 (2020), 266-271.
\arXiv 1904.07350

[Mi12a]
I. Mineyev,
Submultiplicativity and the Hanna Neumann conjecture,
Ann. Math., 175 (2012), 393-414.

[Mi12b]
I. Mineyev,
Groups, graphs, and the Hanna Neumann conjecture,
J. Topol. Anal., 4:1 (2012), 1-12.

[Nos16]
G. A. Noskov,
Mineyev--Dicks proof of the HN-conjecture and the
Euler-Poincar\'e characteristic,
Math. Notes, 99:3 (2016), 390-396.

[Vi49]
A. A. Vinogradov,
On the free products of ordered groups,
Mat. Sb., 25(67):1 (1949), 163-168.

[Za14]
A. Zakharov,
On the rank of the intersection of free subgroups in virtually free groups,
J. Algebra, 418 (2014), 29-43.
\arXiv 1301.3115

\end